\documentclass[12pt]{article}
\usepackage[centertags]{amsmath}
\usepackage{amsfonts}
\usepackage{amssymb}
\usepackage{amsthm}
\usepackage{newlfont}

\date{}

\newcommand{\Z}{\Bbb Z}

\def\NL{\hfill\break}

\begin{document}

\centerline{}

\centerline{}

\centerline {\Large{\bf A note on $\sigma$-reversibility and $\sigma$-symmetry}}

\centerline{}

\centerline{\Large{\bf of skew power series rings}}

\centerline{}
\centerline{}

\centerline{}

\centerline{\bf {L'moufadal Ben Yakoub and Mohamed Louzari}}

\centerline{}

\centerline{Department of mathematics}

\centerline{Abdelmalek Essaadi University}

\centerline{B.P. 2121 Tetouan, Morocco}

\centerline{benyakoub@hotmail.com,  mlouzari@yahoo.com}

\newtheorem{Theorem}{\quad Theorem}[section]

\newtheorem{Definition}[Theorem]{\quad Definition}

\newtheorem{Corollary}[Theorem]{\quad Corollary}

\newtheorem{Lemma}[Theorem]{\quad Lemma}

\newtheorem{Example}[Theorem]{\quad Example}

\newtheorem{Proposition}[Theorem]{\quad Proposition}

\begin{abstract}Let $R$ be a ring and $\sigma$ an endomorphism of $R$. In this note, we study the transfert of the symmetry ($\sigma$-symmetry) and reversibility ($\sigma$-reversibility) from $R$ to its skew power series ring $R[[x;\sigma]]$. Moreover, we study on the relationship between the Baerness, quasi-Baerness and p.p.-property of a ring $R$ and these of the skew power series ring $R[[x;\sigma]]$ in case $R$ is right $\sigma$-reversible. As a consequence we obtain a generalization of \cite{hong/2000}.
\end{abstract}

{\bf Mathematics Subject Classification:} 16S36; 16W20; 16U80\\

{\bf Keywords:} Armendariz rings; Baer rings; p.p.-rings; quasi-Baer rings; skew power series rings; reversible rings; symmetric rings \\

\section{Introduction}

Throughout this paper $R$ denotes an associative ring with identity and $\sigma$ denotes a nonzero non identity endomorphism of a given ring.
\par
\smallskip
Recall that a ring is {\it reduced} if it has no nonzero nilpotent elements. Lambek \cite{lambek}, called a ring $R$ {\it symmetric} if $abc=0$ implies $acb=0$ for $a,b,c\in R$. Every reduced ring is symmetric (\cite[Lemma 1.1]{shin}) but the converse does not hold by \cite[Example II.5]{anderson}. Cohen \cite{cohen}, called a ring $R$ {\it reversible} if $ab=0$ implies $ba=0$ for $a,b\in R$. It is obvious that commutative rings are symmetric and symmetric rings are reversible, but the converse does not hold by \cite[Examples I.5 and II.5]{anderson} and \cite[Examples 5 and 7]{marks}. From \cite{baser/2007}, a ring $R$ is called {\it right $($left$)$ $\sigma$-reversible} if whenever $ab=$ for $a,b\in R$, $b\sigma(a)=0$ ($\sigma(b)a=0$). $R$ is called $\sigma$-reversible if it is both right and left $\sigma$-reversible. Also, by \cite{kwak}, a ring $R$ is called {\it right $($left$)$ $\sigma$-symmetric} if whenever $abc=0$ for $a, b, c \in R$, $ac\sigma(b)=0$ ($\sigma(b)ac=0$). $R$ is called $\sigma$-symmetric if it is both right and left $\sigma$-symmetric. Clearly right $\sigma$-symmetric rings are right $\sigma$-reversible.
\par
\smallskip
Rege and Chhawchharia \cite{chhawchharia}, called a ring $R$ an {\it Armendariz} if whenever polynomials $f=\sum_{i=0}^{n}a_ix^i,\;g=\sum_{j=0}^{m}b_jx^j\in R[x]$ satisfy $fg=0$, then $a_ib_j=0$ for each $i,j$. The Armendariz property of a ring was extended to one of skew polynomial ring in \cite{hong/2003}. For an endomorphism $\sigma$ of a ring $R$, a {\it skew polynomial} ring (also called an {\it Ore extension of endomorphism type}) $R[x;\sigma]$ of $R$ is the ring obtained by giving the polynomial ring over $R$ with the new multiplication $xr=\sigma(r)x$ for all $r\in R$. Also, a {\it skew power series ring} $R[[x;\sigma]]$ is the ring consisting of all power series of the form $\sum_{i=0}^{\infty}a_ix^i\;(a_i\in R)$, which are multiplied using the distributive law and the Ore commutation rule $xa=\sigma(a)x$, for all $a\in R$. According to Hong et al. \cite{hong/2003}, a ring $R$ is called $\sigma$-skew Armendariz if whenever polynomials $f=\sum_{i=0}^{n}a_ix^i$ and $g=\sum_{j=0}^{m}b_jx^j$ $\in R[x;\sigma]$ satisfy $fg=0$ then $a_i\sigma^i(b_j)=0$ for each $i,j$. Baser et al. \cite{baser/2008}, introduced the concept of $\sigma$-(sps) Armendariz rings. A ring $R$ is called $\sigma$-(sps) {\it Armendariz} if whenever $pq=0$ for $p=\sum_{i=0}^{\infty}a_ix^i,\;q=\sum_{j=0}^{\infty}b_jx^j\in R[[x;\sigma]]$, then $a_ib_j=0$ for all $i$ and $j$. According to Krempa \cite{krempa}, an endomorphism $\sigma$ of a ring $R$ is called {\it rigid} if $a\sigma(a)=0$ implies $a=0$ for all $a\in R$. We call a ring $R$ $\sigma$-{\it rigid} if there exists a rigid endomorphism $\sigma$ of $R$. Note that any rigid endomorphism of a ring $R$ is a monomorphism and $\sigma$-rigid rings are reduced by Hong et al. \cite{hong/2000}. Also, by \cite[Theorem 2.8(1)]{kwak}, a ring $R$ is $\sigma$-rigid if and only if $R$ is semiprime right $\sigma$-symmetric and $\sigma$ is a monomorphisme, so right $\sigma$-symmetric ($\sigma$-reversible) rings are a generalization of  $\sigma$-rigid rings.
\par
\smallskip
In this note, we introduce the notion of {\it $\sigma$-skew $($sps$)$ Armendariz} rings which is a generalization of $\sigma$-(sps) Armendariz rings, and we study the transfert of the symmetry ($\sigma$-symmetry) and reversibility ($\sigma$-reversibility) from $R$ to its skew power series ring $R[[x;\sigma]]$. Also we show that $R$ is $\sigma$-(sps) Armendariz if and only if $R$ is $\sigma$-skew (sps) Armendariz and $a\sigma(b)=0$ implies $ab=0$ for $a,b\in R$. Moreover, we study on the relationship between the Baerness, quasi-Baerness and p.p.-property of a ring $R$ and these of the skew power series ring $R[[x;\sigma]]$ in case $R$ is right $\sigma$-reversible. As a consequence we obtain a generalization of \cite{hong/2000}.

\section{$\sigma$-Reversibility and $\sigma$-Symmetry of Skew Power Series Rings}

We introduce the next definition.

\begin{Definition}Let $R$ be a ring and $\sigma$ an endomorphism of $R$. A ring $R$ is called $\sigma$-skew $($sps$)$ Armendariz if whenever $pq=0$ for $p=\sum_{i=0}^{\infty}a_ix^i,\;q=\sum_{j=0}^{\infty}b_jx^j\in R[[x;\sigma]]$, then $a_i\sigma^i(b_j)=0$ for all $i$ and $j$.
\end{Definition}

Every subring $S$ with $\sigma(S)\subseteq S$ of an $\sigma$-skew (sps) Armendariz ring is a $\sigma$-skew (sps) Armendariz ring. In the next, we give an example of a ring $R$ which is $\sigma$-skew (sps) Armendariz but not $\sigma$-(sps) Armendariz.

\begin{Example}\label{ex2}Let $R$ be the polynomial ring $\Z_2[x]$ over $\Z_2$, and let the endomorphism $\sigma\colon R\rightarrow R$ be defined by $\sigma(f(x))=f(0)$ for $f(x)\in \Z_2[x]$.
\NL$(i)$ $R$ is not $\sigma$-$(sps)$ Armendariz because $\sigma$ is not a monomorphism.
\NL$(ii)$ $R$ is an $\sigma$-skew $($sps$)$ Armendariz ring $($as in \cite[Example 5]{hong/2003}$)$. Consider $R[[y;\sigma]]=\Z_2[x][[y;\sigma]]$. Let $p=\sum_{i=0}^{\infty}f_iy^i$ and $q=\sum_{j=0}^{\infty}g_jy^j\in R[[y;\sigma]]$. We have $pq=\sum_{\ell\geq 0}\sum_{\ell=i+j}f_i\sigma^i(g_j)y^{\ell}=0$. If $pq=0$ then $\sum_{\ell=i+j}f_i\sigma^i(g_j)y^{\ell}=0$, for each $\ell\geq 0$. Suppose that there is $f_s\neq 0$ for some $s\geq 0$ and $f_0=f_1=\cdots=f_{s-1}=0$, then $\sum_{i+j=s}f_i\sigma^i(g_j)y^{i+j}=0\Rightarrow f_s\sigma^s(g_0)=0$, since $R$ is a domain then $g_0(0)=0$. Also $\sum_{i+j=s+1}f_i\sigma^i(g_j)y^{i+j}=0\Rightarrow f_s\sigma^s(g_1)+f_{s+1}\sigma^{s+1}(g_0)=0$, since $g_0(0)=0$ then $f_s\sigma^s(g_1)=0$ and so $g_1(0)=0$ by the same method as above. Continuing this process, we have $g_j(0)=0$ for all $j\geq 0$. Thus $f_i\sigma^i(g_j)=0$ for all $i,j$.
\end{Example}

We say that $R$ satisfies the condition $(\mathcal{C_{\sigma}})$, if whenever $a\sigma(b)=0$ for $a,b\in R$, then $ab=0$. By \cite[Theorem 3.3(3iii)]{baser/2008}, if $R$ is $\sigma$-(sps) Armendariz then it satisfies $(\mathcal{C_{\sigma}})$ (so $\sigma$ is a monomorphism). If $R$ is an $\sigma$-skew (sps) Armendariz ring satisfying the condition $(\mathcal{C_{\sigma}})$ then $R$ is $\sigma$-(sps) Armendariz.

\begin{Theorem}A ring $R$ is $\sigma$-$(sps)$ Armendariz ring if and only if it is $\sigma$-skew $(sps)$ Armendariz and satisfies the condition $(\mathcal{C_{\sigma}})$.
\end{Theorem}

\begin{proof}$(\Leftarrow)$. It is clear. $(\Rightarrow)$. If $R$ is $\sigma$-$(sps)$ Armendariz then it satisfies the condition $(\mathcal{C_{\sigma}})$. It suffices to show that if $R$ is $\sigma$-$(sps)$ Armendariz then it is $\sigma$-skew (sps) Armendariz. The proof is similar as of \cite[Theorem 1.8]{hong/2006}. Let $p=\sum_{i=0}^{\infty}a_ix^i$ and $q=\sum_{j=0}^{\infty}b_jx^j\in R[[x;\sigma]]$ with $pq=0$. Note that $a_jb_j=0$ for all $i$ and $j$. We claim that $a_i\sigma^i(b_j)=0$ for all $i$ and $j$. We have $(a_0+a_1x+\cdots)(b_0+b_1x+\cdots)=0$, then $a_0(b_0+b_1x+\cdots)+(a_1x+a_2x^2\cdots)(b_0+b_1x+\cdots)=0$. Since $a_0b_j=0$ for all $j$, we get
$$\;0=(a_1x+a_2x^2+\cdots)(b_0+b_1x+\cdots)$$
$$0=(a_1+a_2x+\cdots)x(b_0+b_1x+\cdots)$$
$$\qquad\quad\; 0=(a_1+a_2x+\cdots)(\sigma(b_0)x+\sigma(b_1)x^2+\cdots).$$
Put $p_1=a_1+a_2x+\cdots$ and $q_1=\sigma(b_0)x+\sigma(b_1)x^2+\cdots$. Since $p_1q_1=0$ then $a_i\sigma(b_j)=0$ for all $i\geq 1$ and $j\geq 0$. We have, also $$0=a_1(\sigma(b_0)x+\sigma(b_1)x^2+\cdots)+(a_2x+a_3x^2+\cdots)(\sigma(b_0)x+\sigma(b_1)x^2+\cdots).$$ Since $a_1\sigma(b_j)=0$ for all $j$, then
$$0=(a_2x+a_3x^2+\cdots)(\sigma(b_0)x+\sigma(b_1)x^2+\cdots)$$
$$\;\;0=(a_2+a_3x+\cdots)(\sigma^2(b_0)x^2+\sigma^2(b_1)x^3+\cdots).$$
Put $p_2=a_2+a_3x+a_4x ^2+\cdots$ and $q_2=\sigma^2(b_0)x^2+\sigma^2(b_1)x^3+\cdots$, and then $p_2q_2=0$ implies $a_i\sigma^2(b_j)=0$ for all $i\geq 2$ and $j\geq 0$. Continuing this process, we can show that $a_i\sigma^i(b_j)=0$ for all $i\geq 0$ and $j\geq 0$. Thus $R$ is $\sigma$-skew $(sps)$ Armendariz.
\end{proof}

\begin{Lemma}\label{lem1}Let $R$ be an $\sigma$-$(sps)$ Armendariz ring. Then for $f=\sum_{i=0}^{\infty}a_ix^i$, $g=\sum_{j=0}^{\infty}b_jx^j$ and $h=\sum_{k=0}^{\infty}c_kx^k\in R[[x;\sigma]]$, if $fgh=0$ then $a_ib_jc_k=0$ for all $i,j,k$.
\end{Lemma}

\begin{proof}Note that, if $fg=0$ then $a_ig=0$ for all $i$. Suppose that $fgh=0$ then $a_i(gh)=0$ for all $i$, and so $(a_ig)h=0$ for all $i$. Therefore $a_ib_jc_k=0$ for all $i,j,k$.
\end{proof}

\begin{Proposition}\label{prop1}Let $R$ be an $\sigma$-$(sps)$ Armendariz ring. Then
\par$(1)$ $R$ is reversible if and only if $R[[x;\sigma]]$ is reversible.
\par$(2)$ $R$ is symmetric if and only if $R[[x;\sigma]]$ is symmetric.
\end{Proposition}

\begin{proof}If $R[[x;\sigma]]$ is symmetric (reversible) then $R$ is symmetric (reversible). Conversely, $(1)$. Let $f=\sum_{i=0}^{\infty}a_ix^i$ and $g=\sum_{j=0}^{\infty}b_jx^j\in R[[x;\sigma]]$, if $fg=0$ then $a_ib_j=0$ for all $i$ and $j$. By \cite[Theorem 3.3 (3ii)]{baser/2008}, we have $\sigma^j(a_i)b_j=0$ for all $i$ and $j$. Since $R$ is reversible, we obtain  $b_j\sigma^j(a_i)=0$ for all $i$ and $j$. Thus $gf=\sum_{\ell=0}^{\infty}\sum_{\ell=i+j}b_j\sigma^j(a_i)x^{\ell}=0$. $(2)$. We will use freely \cite[Theorem 3.3 (3ii)]{baser/2008}, reversibility and symmetry of $R$. Let $f=\sum_{i=0}^{\infty}a_ix^i$, $g=\sum_{j=0}^{\infty}b_jx^j$ and $h=\sum_{k=0}^{\infty}c_kx^k\in R[[x;\sigma]]$, if $fgh=0$ then $a_ib_jc_k=0$ for all $i$, $j$ and $k$, by Lemma \ref{lem1}. Then for all $i,j,k$ we have $b_jc_ka_i=0\Rightarrow \sigma^k(b_j)c_ka_i=0 \Rightarrow a_i\sigma^k(b_j)c_k=0\Rightarrow a_ic_k\sigma^k(b_j)=0\Rightarrow c_k\sigma^k(b_j)a_i=0\Rightarrow \sigma^i[c_k\sigma^k(b_j)]a_i=0\Rightarrow a_i\sigma^i[c_k\sigma^k(b_j)]=0$. Thus $fhg=0$.
\end{proof}

The next Lemma gives a relationship between $\sigma$-reversibility ($\sigma$-symmetry) and reversibility (symmetry).

\begin{Lemma}[{\cite[Lemma 3.1]{louzari2}}]\label{lem2}Let $R$ be a ring and $\sigma$ an endomorphism of $R$. If $R$ satisfies the condition $(\mathcal{C_{\sigma}})$. Then
\par $(1)$ $R$ is reversible if and only if $R$ is $\sigma$-reversible;
\par $(2)$ $R$ is symmetric if and only if $R$ is $\sigma$-symmetric.
\end{Lemma}

\begin{Theorem}\label{s rev theo}Let $R$ be an $\sigma$-$(sps)$ Armendariz ring. The following statements are equivalent:
\par$(1)$ $R$ is reversible $($symmetric$)$;
\par$(2)$ $R$ is $\sigma$-reversible $($$\sigma$-symmetric$)$;
\par$(3)$ $R$ is right $\sigma$-reversible $($right $\sigma$-symmetric$)$;
\par$(4)$ $R[[x;\sigma]]$ is reversible $($symmetric$)$.
\end{Theorem}

\begin{proof}We prove the reversible case (the same for the symmetric case).\NL$(1)\Leftrightarrow (4)$. By Proposition \ref{prop1}.
\NL$(1)\Rightarrow (2)$ and $(2)\Rightarrow (3)$. Immediately from Lemma \ref{lem2}.
\NL$(3)\Rightarrow (1)$. Let $a,b\in R$, if $ab=0$ then $b\sigma(a)=0$ (right $\sigma$-reversibility), so $ba=0$ (condition $(\mathcal{C_{\sigma}})$).
\end{proof}

\section{Related Topics}

In this section we turn our attention to the relationship between the Baerness, quasi-Baerness and p.p.-property of a ring $R$ and these of the skew power series ring $R[[x;\sigma]]$ in case $R$ is right $\sigma$-reversible. For a nonempty subset $X$ of $R$, we write $r_R(X)=\{c\in R|dc=0\;\mathrm{for\;any}\;d\in X\}$ which is called the right annihilator of $X$ in $R$.

\begin{Lemma}\label{lem4}If $R$ is a right $\sigma$-reversible ring with $\sigma(1)=1$. Then
\par$(1)$ $\sigma(e)=e$ for all idempotent $e\in R$;
\par$(2)$ $R$ is abelian.
\end{Lemma}

\begin{proof}(1) Let $e$ an idempotent of $R$. We have $e(1-e)=(1-e)e=0$ then $(1-e)\sigma(e)=e\sigma((1-e))=0$, so $\sigma(e)-e\sigma(e)=e-e\sigma(e)=0$, therefore $\sigma(e)=e$. (2) Let $r\in R$ and $e$ an idempotent of $R$. We have $e(1-e)=0$ then $e(1-e)r=0$, since $R$ is right $\sigma$-reversible then $(1-e)r\sigma(e)=0=(1-e)re=0$, so $re=ere$. Since $(1-e)e=0$, we have also $er=ere$. Then $R$ is abelian.
\end{proof}

\begin{Lemma}\label{lem3}Let $R$ be a right $\sigma$-reversible ring with $\sigma(1)=1$, then the set of all idempotents in $R[[x;\sigma]]$ coincides with the set of all idempotents of $R$. In this case $R[[x;\sigma]]$ is abelian.
\end{Lemma}

\begin{proof}We adapt the proof of \cite[Theorem 2.13(iii)]{baser/2007} for $R[[x;\sigma]]$. Let $f^2=f\in R[[x;\sigma]]$, where $f=f_0+f_1x+f_2x^2+\cdots $. Then $$\sum_{\ell=0}^{\infty}\sum_{\ell=i+j}f_i\sigma^i(f_j)x^{\ell}=\sum_{\ell=0}^{\infty}f_{\ell}x^{\ell}.$$
For $\ell=0$, we have $f_0^2=f_0$. For $\ell=1$, we have $f_0f_1+f_1\sigma(f_0)=f_1$, but $f_0$ is central and $\sigma(f_0)=f_0$, so $f_0f_1+f_1f_0=f_1$, a multiplication by $(1-f_0)$ on the left hand gives $f_1=f_0f_1$, and so $f_1=0$. For $\ell=2$, we have $f_0f_2+f_1\sigma(f_1)+f_2\sigma^2(f_0)=f_2$, so $f_0f_2+f_2f_0=f_2$ (because $f_1=0$ and $\sigma^2(f_0)=f_0$), a multiplication by $(1-f_0)$ on the left hand gives $f_0f_2=f_2=0$. Continuing this procedure yields $f_i=0$ for all $i\geq 1$. Consequently, $f=f_0=f_0^2\in R$. Since $R$ is abelian then  $R[[x;\sigma]]$ is abelian.
\end{proof}

Kaplansky \cite{Kaplansky}, introduced the concept of {\it Baer rings} as rings in which the right (left) annihilator of every nonempty subset is generated by an idempotent. According to Clark \cite{clark}, a ring $R$ is called {\it quasi-Baer} if the right annihilator of each right ideal of $R$ is generated (as a right ideal) by an idempotent. It is well-known that these two concepts are left-right symmetric. A ring $R$ is called a {\it right} ({\it left}) {\it p.p.-ring} if the right (left) annihilator of an element of $R$ is generated by an idempotent. $R$ is called a {\it p.p.-ring} if it is both a right and left p.p.-ring.

\begin{Theorem}\label{theo baer qbaer}Let $R$ be a right $\sigma$-reversible ring with $\sigma(1)=1$. Then
\par$(1)$ $R$ is a Baer ring  if and only if $R[[x;\sigma]]$ is a Baer ring;
\par$(2)$ $R$ is a quasi-Baer ring  if and only if $R[[x;\sigma]]$ is a quasi-Baer ring.
\end{Theorem}

\begin{proof}$(\Rightarrow)$. Suppose that $R$ is Baer. Let $A$ be a nonempty subset of $R[[x;\sigma]]$ and $A^*$ be the set of all coefficients of elements of $A$. Then $A^*$ is a nonempty subset of $R$ and so $r_R(A^*)=eR$ for some idempotent element $e\in R$. Since $e\in r_{R[[x;\sigma]]}(A)$ by Lemma \ref{lem4}. We have $eR[[x;\sigma]]\subseteq r_{R[[x;\sigma]]}(A)$. Now, let $0\neq q=b_0+b_1x+b_2x^2+\cdots\in r_{R[[x;\sigma]]}(A)$. Then $Aq=0$ and hence $pq=0$ for any $p\in A$. Let $p=a_0+a_1x+a_2x^2+\cdots$, then
$$pq=\sum_{\ell\geq 0}\sum_{\ell=i+j}a_i\sigma^i(b_j)x^{\ell}=0.$$
\begin{itemize}
\item $\ell=0$ implies $a_0b_0=0$ then $b_0\in r_R(A^*)=eR$.
\item $\ell=1$ implies $a_0b_1+a_1\sigma(b_0)=0$, since $b_0=eb_0$ and $\sigma(e)=e$ then $a_0b_1+a_1e\sigma(b_0)=0$, but $a_1e=0$ so $a_0b_1=0$ and hence $b_1\in r_R(A^*)$.
\item $\ell=2$ implies $a_0b_2+a_1\sigma(b_1)+a_2\sigma^2(b_0)=0$, then $a_0b_2+a_1e\sigma(b_1)+a_2e\sigma^2(b_0)=0$, but $a_1e\sigma(b_1)=a_2e\sigma^2(b_0)=0$, hence $a_0b_2=0$. Then $b_2\in r_R(A^*)$.
\end{itemize}
Continuing this procedure yields $b_0,b_1,b_2,b_3,\cdots\in r_R(A^*)$. So, we can write $q=eb_0+eb_1x+eb_2x^2+\cdots\in eR[[x;\sigma]]$. Therefore $eR[[x;\sigma]]=r_{R[[x;\sigma]]}(A)$. Consequently, $R[[x;\sigma]]$ is a Baer ring.\par Conversely, Suppose that $R[[x;\sigma]]$ is Baer. Let $B$ be a nonempty subset of $R$. Then $r_{R[[x;\sigma]]}(B)=eR[[x;\sigma]]$ for some idempotent $e\in R$ by Lemma \ref{lem3}. Thus $r_R(B)=r_{R[[x;\sigma]]}(B)\cap R=eR[[x;\sigma]]\cap R=eR$. Therefore $R$ is Baer.\par The proof for the case of the quasi-Baer property follows in a similar fashion; In fact, for any right ideal $A$ of $R[[x;\sigma]]$, take $A^*$ as the right ideal generated by all coefficients of elements of $A$.
\end{proof}

From \cite[Example 20]{hong/2000}, $R=M_2(\Z)$ is a Baer ring and $R[[x]]$ is not Baer. Clearly $R$ is not reversible. So that, the ``right $\sigma$-reversibility'' condition in Theorem \ref{theo baer qbaer}(1) is not superfluous.

According to Annin \cite{annin}, a ring $R$ is $\sigma$-{\it compatible} if for each $a,b\in R$, $a\sigma(b)=0$ if and only if $ab=0$. Hashemi and Moussavi \cite[Corollary 2.14]{hashemi/quasi} have proved Theorem \ref{theo baer qbaer}(2), when $R$ is $\sigma$-compatible. Consider $R$ and $\sigma$ as in Example \ref{ex2}. Since $R$ is a domain then it is right $\sigma$-reversible (with $\sigma(1)=1$). Also $R$ is not $\sigma$-compatible (so $R$ does not satisfy the condition $(\mathcal{C_{\sigma}})$), because $\sigma$ is not a monomorphism. Therefore Theorem \ref{theo baer qbaer}(2) is not a consequence of \cite[Corollary 2.14]{hashemi/quasi}. On other hand, if $R$ is reversible then $\sigma$-compatibility implies right $\sigma$-reversibility. But, if $R$ is not reversible, we can easily see that this implication does not hold.

\begin{Theorem}\label{theo pp.rings}Let $R$ be a right $\sigma$-reversible ring with $\sigma(1)=1$.  If $R[[x;\sigma]]$ is a p.p.-ring then $R$ is a p.p.-ring.
\end{Theorem}

\begin{proof}Suppose that $R[[x;\sigma]]$ is a right p.p.-ring. Let $a\in R$, then there exists an idempotent $e\in R$ such that $r_{R[[x;\sigma]]}(a)=eR[[x;\sigma]]$ by Lemma \ref{lem3}. Hence $r_R(a)=eR$, and therefore $R$ is a right p.p.-ring.
\end{proof}

Also, in Example \ref{ex2}, $R$ is not $\sigma$-(sps) Armendariz. So Theorem \ref{theo baer qbaer} and Theorem \ref{theo pp.rings} are not consequences of \cite[Theorem 3.2]{baser/2008}.

Since $\sigma$-rigid rings are right $\sigma$-reversible \cite[Theorem 2.8 (1)]{kwak}, we have the following Corollaries.

\begin{Corollary}[{\cite[Theorem 21]{hong/2000}}]Let $R$ be an $\sigma$-rigid ring. Then $R$ is a Baer ring if and only if $R[[x;\sigma]]$ is a Baer ring.
\end{Corollary}

\begin{Corollary}[{\cite[Corollary 22]{hong/2000}}]Let $R$ be an $\sigma$-rigid ring. Then $R$ is a quasi-Baer ring if and only if $R[[x;\sigma]]$ is a quasi-Baer ring.
\end{Corollary}

{\bf ACKNOWLEDGEMENTS.} The second author wishes to thank Professor Amin Kaidi of University of Almer\'ia for his generous hospitality. This work was supported by the project PCI Moroccan-Spanish A/011421/07.

\end{document}